\documentclass[11pt]{article}
\usepackage[a4paper,margin=1in]{geometry}
\usepackage[T1]{fontenc}
\usepackage[utf8]{inputenc}
\usepackage{amsmath,amssymb,amsthm,mathtools}
\usepackage{enumitem}
\usepackage{microtype}
\usepackage[colorlinks=true,linkcolor=blue,citecolor=blue,urlcolor=blue]{hyperref}

\newtheorem{theorem}{Theorem}[section]
\newtheorem{proposition}[theorem]{Proposition}
\newtheorem{lemma}[theorem]{Lemma}
\newtheorem{remark}[theorem]{Remark}
\newenvironment{keyword}{\par\vspace{0.5em}\noindent\textit{Keywords:}\ }{\par\vspace{0.5em}}

\let\Ent\relax
\DeclareMathOperator{\Ent}{Ent}
\DeclareMathOperator{\GL}{GL}
\DeclareMathOperator{\SL}{SL}
\DeclareMathOperator{\spanop}{span}

\title{The Asymptotically Sharp logarithmic Sobolev Constant for the All-Transvection Walk on $\mathrm{GL}_n(\mathbb F_2)$}
\author{Xumin Wang}
\date{}

\begin{document}
\maketitle

\begin{abstract}
Let $T_n$ be the conjugacy class in $\GL_n(\mathbb F_2)$ consisting of all nontrivial transvections. This paper studies the  logarithmic Sobolev constant of the random walk whose transition distribution is the uniform probability measure on $T_n$. Using column-wise conditional entropy decomposition, the optimal logarithmic Sobolev constant of the loopless complete graph, Hahn--Banach extension, transvection connections between conditional slices, and averaging over left and right translations, we prove that
\[
C_{\mathrm{LS}}(n)=n^2\log 2+O(1).
\]
\end{abstract}


\section{Introduction}
The logarithmic Sobolev (abbreviated as log--Sobolev) inequality is an important tool for studying entropy decay, hypercontractivity, and mixing times of finite Markov chains. Gross's classical work \cite{Gross1975} established the fundamental connection between the log--Sobolev inequality and hypercontractivity. For reversible Markov chains on finite state spaces, Diaconis and Saloff-Coste \cite{DSC1996} systematically developed the log--Sobolev inequality method and determined the optimal constants for several basic models.

There has been a great deal of work on optimal hypercontractivity for groups and group von Neumann algebras. In Fourier analysis in the commutative setting, Bonami \cite{Bonami1970} and Beckner \cite{Beckner1975} established several hypercontractive fundamental inequalities. Andersson's studies of hypercontractivity for finite groups of low order and Poisson integrals on finite groups \cite{Andersson1999,Andersson2002} are closely related to the case of cyclic groups of very small order. 
In the noncommutative setting, Biane \cite{Biane1997} proved that the $q$-Ornstein--Uhlenbeck semigroup on $q$-Gaussian von Neumann algebras has optimal hypercontractivity. Junge, Palazuelos, Parcet, and Perrin developed a method for studying hypercontractivity of group von Neumann algebras, including the cases of free groups and finite cyclic groups \cite{JPPP2017}. Junge, Palazuelos, Parcet, Perrin, and Ricard studied hypercontractivity for free products and, in particular, gave hypercontractivity estimates for the free Poisson semigroup on free group von Neumann algebras \cite{JPPPR2015}. In recent years, There is already extensive results for finite cyclic groups, as follows. Cao, Fan, Han, Qiu, and Wang \cite{CFHQW2026} determined the optimal hypercontractivity constants for $\mathbb Z_3$ and for biased Bernoulli random variables. For $\mathbb Z_n$, $n\ge4$, Xie and Zhang \cite{XieZhang2026} and Yao \cite{Yao2026} independently proved that the optimal log--Sobolev constant for the Poisson semigroup on $\mathbb Z_n$ determined by the word length function $\psi_n(k)=\min(k,n-k)$ is $2$. For the nearest-neighbor simple random walk on the cycle graph, Frank and Ivanisvili \cite{FI2026} proved that, for $n\ge4$, the optimal constant in the form of the inequality used in this paper equals twice the reciprocal of the spectral gap. Xie and Zhang also established quartic stability estimates. Chen, Lam, and Lu \cite{CLL2026} further determined the optimal quartic stability coefficient that holds uniformly for all $n\ge4$.

Random walks on finite groups are closely related to Markov chain theory, representation theory, and asymptotic combinatorics. Diaconis and Shahshahani \cite{DS1981} used representation theory to study the random transposition walk on the symmetric group. This model remains a typical model in the study of random walks on finite groups. For a general theory of random walks on finite groups, see the survey by Saloff-Coste \cite{SaloffCoste2004}. Among linear groups, the random transvection walk is a natural generalization of the random transposition walk: at each step, instead of multiplying by a transposition, one multiplies by a linear transformation whose difference from the identity is a rank-one nilpotent operator. Hildebrand \cite{Hildebrand1992} studied the random transvection walk on $\SL_n(\mathbb F_q)$. In his model, at each step one chooses uniformly a group element from the conjugacy class consisting of all nontrivial transvections. To distinguish this model from the random walk generated by elementary transvections with respect to a fixed coordinate system, we call this model the all-transvection random walk. Hildebrand proved that this random walk exhibits cutoff around time $n$. J.~He \cite{He2020} studied a related random walk in the space of symplectic forms over finite fields and proved the corresponding cutoff result.

With respect to a fixed coordinate system, the random walk generated by elementary transvections is a closely related model with weaker symmetry. On $\GL_n(\mathbb F_2)$, each step of this Markov chain selects two distinct rows (with order distinguished) and adds one row to the other. Ben-Hamou \cite{BenHamou2025} proved a log--Sobolev inequality for this random walk with constant of order $n^2$, and thereby obtained an upper bound $O(n^2\log n)$ for the mixing time. Pillai and Smith \cite{PillaiSmith2026} studied the projection of the elementary transvection random walk onto the first $k$ columns and, combining local-to-global entropy estimates with burn-in argument, proved that the mixing time of the first $k$ columns is $O(nk\log n)$, thus resolving a conjecture of Ben-Hamou and Peres.

This paper studies the log--Sobolev constant of Hildebrand's random transvection model in the case $q=2$. Our goal is to go beyond the existing cutoff theory and obtain precise entropy dissipation estimates, proving that the optimal log--Sobolev constant is asymptotically equal to $\log|G_n|$.

We now define the random walk model studied in this paper. Throughout, assume that $n\ge 2$. Let $\mathbb F_2$ be the field with two elements, let $V=\mathbb F_2^n$, and 
\[
G_n=\GL(V)=\GL_n(\mathbb F_2)=\SL_n(\mathbb F_2),
\]
where the last equality holds because $\mathbb F_2^\times=\{1\}$. Let $V^*$ be the dual space of $V$. For $u\in V$ and $\varphi\in V^*\setminus\{0\}$, define the rank-one operator
\[
(u\otimes \varphi)(x)=\varphi(x)u,\ \ x\in V.
\]
From \cite{Hildebrand1992}, following Hildebrand's random transvection model, define
\[
T_n=\{I+u\otimes \varphi: u\in V\setminus\{0\},\ \varphi\in V^*\setminus\{0\},\ \varphi(u)=0\}.
\]
$T_n$ is the all-transvection conjugacy class in $G_n$. If $\tau=I+u\otimes \varphi\in T_n$, then
\[
(u\otimes \varphi)^2=\varphi(u)u\otimes \varphi=0,\ \ \tau^2=I.
\]
Hence $T_n=T_n^{-1}$ and $I\notin T_n$. Moreover,
\[
|T_n|=(2^n-1)(2^{n-1}-1).
\]

Let $\pi_n$ be the uniform probability measure on $G_n$, and $\mu_n$ be the uniform probability measure on $T_n$, namely,
\[
\pi_n(g)=|G_n|^{-1},\ \ \mu_n(\tau)=|T_n|^{-1}.
\]
And let 
\[
L^2(G_n,\pi_n)
=\left\{ f:G_n\to\mathbb R \;\middle|\; 
\langle f,f\rangle_{\pi_n}<\infty
\right\},
\]
with inner product
\[
\langle f_1,f_2\rangle_{\pi_n}=\sum_{g\in G_n} f_1(g)f_2(g)\pi_n(g).
\]
Define the Markov operator $P_n:L^2(G_n,\pi_n)\to L^2(G_n,\pi_n)$ by
\[
P_nf(g)=\sum_{\tau\in T_n} f(\tau g)\mu_n(\tau).
\]
Since $T_n=T_n^{-1}$ and $\pi_n$ is left-invariant, $P_n$ is a self-adjoint Markov operator on $L^2(G_n,\pi_n)$. Its Dirichlet form is
\[
\mathcal E_n(f,f)=\langle f,(I-P_n)f\rangle_{\pi_n}.
\]
Equivalently, on the product probability space $(T_n\times G_n,\mu_n\otimes \pi_n)$, we also have
\[
\mathcal E_n(f,f)=\frac12\sum_{\tau\in T_n}\sum_{g\in G_n}
\bigl(f(\tau g)-f(g)\bigr)^2\mu_n(\tau)\pi_n(g).
\]

For a nonnegative function $F:G_n\to [0,\infty)$, define
\[
\Ent_{\pi_n}(F)
=\sum_{g\in G_n}F(g)\log F(g)\pi_n(g)
-\left(\sum_{g\in G_n}F(g)\pi_n(g)\right)
\log\left(\sum_{g\in G_n}F(g)\pi_n(g)\right),
\]
with the convention $0\log 0=0$. We recall the log-Sobolev inequality,
\[
\Ent_{\pi_n}(f^2)\le C\mathcal E_n(f,f).
\]
The optimal constant is denoted by
\[
C_{\mathrm{LS}}(n)=\sup_{f\not\equiv \mathrm{const}}
\frac{\Ent_{\pi_n}(f^2)}{\mathcal E_n(f,f)}.
\]

\begin{theorem}\label{t1}
Let $G_n$ and $C_{\mathrm{LS}}(n)$ be as above, as $n\to\infty$, then
\[
\log |G_n|\le C_{\mathrm{LS}}(n)\le \log |G_n|+o(1).
\]
Consequently,
\[
C_{\mathrm{LS}}(n)=\log |G_n|+o(1)=n^2\log 2+O(1).
\]
\end{theorem}

The proof proceeds as follows. View each matrix as an ordered basis consisting of its column vectors, and consider the increasing family of $\sigma$-algebras generated by the first $k$ columns. Hence the entropy on $G_n$ can be decomposed as a sum of conditional entropies corresponding to the individual columns. When the first $k$ linearly independent column vectors are fixed, the next column is uniformly distributed on the complement of the subspace spanned by the first $k$ linearly independent columns. Therefore, one can use the optimal log--Sobolev constant of the loopless complete graph to estimate the corresponding conditional entropy. With the first $k$ columns held fixed, for any two possible values of the next column there exists a transvection mapping one value to the other. Thus the Dirichlet form on the complete graph is controlled by the local energy associated with the corresponding transvection. Finally, by averaging over left and right translations, one shows that the averaged local energy is exactly equal to the original Dirichlet form $\mathcal E_n$.

\section{Column-wise entropy decomposition and local complete-graph estimates}

Every $g\in G_n$ can be written uniquely as an ordered basis consisting of column vectors,
\[
g=(v_1,\ldots,v_n),\ \ v_i\in V.
\]
This section proceeds in two steps. First, we record the exact column counting compute the corresponding martingale conditional expectations over completion sets. Then, using the entropy estimate on the complete graph, the action of transvections between different choices of columns, and the reverse triangle inequality for the $L^2$ norm, we prove the local estimate needed later.

\subsection{Counting columns and averages over completion sets}\label{s21}

Counting the number of ordered bases column by column, we obtain
\[
|G_n|=\prod_{k=1}^n m_k,\ \ m_k=2^n-2^{k-1}.
\]
Equivalently,
\[
|G_n|=\prod_{j=0}^{n-1}(2^n-2^j)
=2^{n^2}\prod_{r=1}^n(1-2^{-r}),
\]
and hence
\[
\log |G_n|=n^2\log 2+\sum_{r=1}^n\log(1-2^{-r})
=n^2\log 2+O(1).
\]

On the probability space $(G_n,\pi_n)$, let $V_k:G_n\to V$, namely, $V_k(g)=v_k$. For $0\le k\le n$, define
\[
\mathcal F_k=\sigma(V_1,\ldots,V_k),\ \ \mathcal F_0=\{\emptyset,G_n\}.
\] 
If $v_1,\ldots,v_{k-1}$ are linearly independent, let  $v_{<k}=(v_1,\ldots,v_{k-1})$ and
\[
U(v_{<k})=\spanop(v_1,\ldots,v_{k-1}),\ \ X(v_{<k})=V\setminus U(v_{<k}).
\]
Then
\[
|X(v_{<k})|=m_k=2^n-2^{k-1}.
\]
Therefore, the first $k-1$ columns are fixed to be $v_{<k}$, the random column $V_k$ is uniformly distributed over $X(v_{<k})$. For every $a\in X(v_{<k})$, the number of ordered bases whose first $k$ columns are $(v_1,\ldots,v_{k-1},a)$ is
\[
\prod_{j=k+1}^n (2^n-2^{j-1}),
\]
with the empty product interpreted as $1$. Obviously, this number is independent of $a$.

The above counting also describes the completions once the $k$-th column is fixed. For $a\in X(v_{<k})$, let
\[
B(v_{<k},a)=\{(v_1,\ldots,v_{k-1},a,v_{k+1},\ldots,v_n)\in G_n\}.
\]
Then
\[
|B(v_{<k},a)|=\prod_{j=k+1}^n(2^n-2^{j-1})
=\prod_{j=k+1}^n m_j,
\]
with the usual convention that the empty product equals $1$. Hence, 
\begin{align*}
	\mathbb P_{\pi_n}(V_k=a\mid V_{<k}=v_{<k})
	&=\frac{|B(v_{<k},a)|}{\sum_{b\in X(v_{<k})}|B(v_{<k},b)|}=\frac1{m_k},\\
	\mathbb P_{\pi_n}(g=h\mid V_{<k}=v_{<k},V_k=a)
	&=\frac1{|B(v_{<k},a)|},\ \  h\in B(v_{<k},a).
\end{align*}

Let $F=f^2$, and define the Doob martingale
\[
M_k=\mathbb E_{\pi_n}[F\mid \mathcal F_k], \ \ 0\le k\le n.
\]
We have $M_0=\mathbb E_{\pi_n}F$ and $M_n=F$. By the uniformity above, we obtain an explicit formula for $M_k$. If $g\in B(v_{<k},a)$, then
\[
M_k(g)
=\frac{1}{|B(v_{<k},a)|}
\sum_{h\in B(v_{<k},a)} f(h)^2.
\]
Denote the square root of this average by
\[
\Psi_{v_{<k}}(a)
=
\left(
\frac{1}{|B(v_{<k},a)|}
\sum_{h\in B(v_{<k},a)} f(h)^2
\right)^{1/2}.
\]
Thus, on the set $B(v_{<k},a)$, we have
\[
M_k=\Psi_{v_{<k}}(a)^2.
\]
Similarly, on the set of matrices determined by the first $k-1$ columns $v_{<k}$, we have
\[
M_{k-1}
=\frac1{m_k}\sum_{a\in X(v_{<k})}\Psi_{v_{<k}}(a)^2.
\]

For a nonnegative random variable $Y$ and a sub-$\sigma$-algebra $\mathcal G$, define the conditional entropy by
\[
\Ent(Y\mid \mathcal G)
=\mathbb E[Y\log Y\mid \mathcal G]
-\mathbb E[Y\mid \mathcal G]\log \mathbb E[Y\mid \mathcal G].
\]

\begin{proposition}\label{p21}
For every nonnegative $F$ on $G_n$, if
\[
M_k=\mathbb E_{\pi_n}[F\mid \mathcal F_k],\ \ 0\le k\le n,
\]
then
\[
\Ent_{\pi_n}(F)
=\sum_{k=1}^n \mathbb E_{\pi_n}\bigl[\Ent(M_k\mid \mathcal F_{k-1})\bigr].
\]
\end{proposition}

\begin{proof}
Since $M_0=\mathbb E_{\pi_n}F$ and $M_n=F$, and from telescoping, then
\begin{align*}
\Ent_{\pi_n}(F)
&=\mathbb E_{\pi_n}[M_n\log M_n]-M_0\log M_0 \\
&=\sum_{k=1}^n
\mathbb E_{\pi_n}\bigl[M_k\log M_k-M_{k-1}\log M_{k-1}\bigr].
\end{align*}
For each $k$, we have $M_{k-1}=\mathbb E_{\pi_n}[M_k\mid \mathcal F_{k-1}]$. Hence the $k$-th summand is
\[
\mathbb E_{\pi_n}\bigl[
\mathbb E[M_k\log M_k\mid \mathcal F_{k-1}]
-\mathbb E[M_k\mid \mathcal F_{k-1}]\log \mathbb E[M_k\mid \mathcal F_{k-1}]
\bigr],
\]
which is exactly $\mathbb E_{\pi_n}[\Ent(M_k\mid \mathcal F_{k-1})]$. Summing over $k$ proves the proposition.
\end{proof}

\subsection{From complete-graph entropy to selected all-transvection energies}\label{s22}

We first introduce the all-transvection energy that will be used to estimate the column energy on the complete graph. For each $k$, assume for the moment that a map $\Theta_k$ has been chosen satisfying the following condition: for every triple $(v_{<k},a,b)$ with $a,b\in X(v_{<k})$ and $a\ne b$, we have
\[
\Theta_k(v_{<k},a,b)\in T_n
\]
and
\[
\Theta_k(v_{<k},a,b)|_{U(v_{<k})}=\mathrm{id},\ \ \Theta_k(v_{<k},a,b)a=b.
\]
Lemma~\ref{l24} will show that such a choice is always possible.

Once $\Theta_k$ is fixed, define
\[
\Omega_k=
\{(g,b): g=(v_1,\ldots,v_n)\in G_n,\ b\in X(v_{<k})\setminus\{v_k\}\}.
\]
Let $\rho_k$ be the probability measure on $\Omega_k$ determined by
\[
\sum_{(g,b)\in\Omega_k}\Phi(g,b)\rho_k(g,b)
=
\sum_{g=(v_1,\ldots,v_n)\in G_n}\pi_n(g)
\frac1{m_k-1}
\sum_{b\in X(v_{<k})\setminus\{v_k\}}
\Phi(g,b).
\]
For $(g,b)\in\Omega_k$, let
\[
\tau_k(g,b)=\Theta_k(v_{<k},v_k,b),
\]
and define
\[
Q_k(f)
=\frac12
\sum_{(g,b)\in\Omega_k}
\bigl(f(\tau_k(g,b)g)-f(g)\bigr)^2\rho_k(g,b).
\]
The goal of this section is to prove the following inequality, which controls the global entropy by local energies,
\begin{equation}\tag{2.1}\label{e21}
\Ent_{\pi_n}(f^2)
\le
\sum_{k=1}^n c_{m_k} Q_k(f),
\end{equation}
where $c_m$ is the optimal log-Sobolev constant of the loopless complete graph on $m$ points. We first express the conditional entropy corresponding to each column as an entropy on a complete graph, and then estimate differences of function values on the complete graph by the differences of function values associated with the chosen transvection in $Q_k(f)$.

Let $\nu_{v_{<k}}$ denote the uniform probability measure on $X(v_{<k})$.

\begin{lemma}\label{l22}
For every $v_{<k}$, we have
\[
\Ent(M_k\mid \mathcal F_{k-1})(v_{<k})
=
\Ent_{\nu_{v_{<k}}}\bigl(\Psi_{v_{<k}}^2\bigr).
\]
\end{lemma}

\begin{proof}
Fix $v_{<k}$, then $\mathbb P_{\pi_n}(V_k=a\mid V_{<k}=v_{<k})=m_k^{-1}$. On the set $B(v_{<k},a)$, we have 
\[
M_k=\Psi_{v_{<k}}(a)^2.
\]
Therefore
\[
\mathbb E[M_k\log M_k\mid V_{<k}=v_{<k}]
=\frac1{m_k}\sum_{a\in X(v_{<k})}
\Psi_{v_{<k}}(a)^2\log \Psi_{v_{<k}}(a)^2,
\]
and
\[
\mathbb E[M_k\mid V_{<k}=v_{<k}]
=\frac1{m_k}\sum_{a\in X(v_{<k})}
\Psi_{v_{<k}}(a)^2.
\]
This value is exactly the entropy of $\Psi_{v_{<k}}^2$ with respect to the uniform probability measure $\nu_{v_{<k}}$.
\end{proof}

Next we estimate the above one-column conditional entropy using the Dirichlet form on the complete graph. Let $X$ be a finite set with $m$ elements, and let $\nu_X$ be the uniform probability measure on $X$. Consider the loopless complete-graph Markov kernel
\[
K_X(x,y)=\frac1{m-1}\mathbf 1_{\{x\ne y\}},\ \ x,y\in X.
\]
The corresponding Dirichlet form is
\[
\mathcal D_X(\phi,\phi)
=\frac12\sum_{x,y\in X}\nu_X(x)K_X(x,y)(\phi(x)-\phi(y))^2
=\frac{1}{2m(m-1)}\sum_{x\ne y}(\phi(x)-\phi(y))^2.
\]
By the result of Diaconis and Saloff-Coste \cite[Appendix A.5]{DSC1996} on the optimal log-Sobolev constant for the complete graph, we have
\begin{equation}\tag{2.2}\label{e22}
	\Ent_{\nu_X}(\phi^2)\le c_m \mathcal D_X(\phi,\phi),\ \
	c_m=\begin{cases}\dfrac{(m-1)\log(m-1)}{m-2},&m>2,\\1,&m=2.\end{cases}
\end{equation}
where $c_2$ is obtained as the limit of the above formula. We shall use the following asymptotic estimate,
\[
c_m=\log m+O\left(\frac{\log m}{m}\right),\ \ m\to\infty.
\]

For $v_{<k}$ and $X(v_{<k})$, define
\[
D_k(v_{<k})
=\frac{1}{2m_k(m_k-1)}
\sum_{\substack{a,b\in X(v_{<k})\\ a\ne b}}
\bigl(\Psi_{v_{<k}}(a)-\Psi_{v_{<k}}(b)\bigr)^2.
\]
From the above inequality on the complete graph and Lemma~\ref{l22}, we obtain the following estimate.
\begin{lemma}\label{l23}
For every real-valued $f$ on $G_n$, we have
\[
\Ent_{\pi_n}(f^2)
\le
\sum_{k=1}^n c_{m_k}\,\mathbb E_{\pi_n}\bigl[D_k(V_{<k})\bigr].
\]
\end{lemma}

\begin{proof}
Since \(|X(v_{<k})|=m_k\) and by \eqref{e22}, then
\[
\Ent_{\nu_{v_{<k}}}(\Psi_{v_{<k}}^2)\le c_{m_k}D_k(v_{<k}).
\]
Combining Proposition~\ref{p21} and Lemma~\ref{l22}, we obtain
\begin{align*}
	\Ent_{\pi_n}(f^2)&=\sum_{k=1}^n\mathbb E_{\pi_n}\bigl[\Ent(M_k\mid\mathcal F_{k-1})\bigr]\\
	&=\sum_{k=1}^n\mathbb E_{\pi_n}\bigl[\Ent_{\nu_{V_{<k}}}(\Psi_{V_{<k}}^2)\bigr]\\
	&\le\sum_{k=1}^n c_{m_k}\mathbb E_{\pi_n}[D_k(V_{<k})].\qedhere
\end{align*}
\end{proof}

\begin{lemma}\label{l24}
Fix  $v_{<k}$, and let $U=U(v_{<k})$. If $a,b\in V\setminus U$ and $a\ne b$, then there exists $\tau\in T_n$ such that
\[
\tau|_U=\mathrm{id}_U,
\qquad
\tau a=b.
\]
\end{lemma}

\begin{proof}
Let $u=a+b\ne0$ and $W=U+\mathbb F_2u$. Since
\[
W=U\cup(u+U),\ \ a\notin U,\ \ a+u=b\notin U,
\]
we have $a\notin W$, and hence $W+\mathbb F_2a=W\oplus\mathbb F_2a$. By the Hahn--Banach extension theorem \cite[Chapter III, Section 6]{Conway1990} for linear functionals, there exists $\varphi\in V^*$ such that
\[
\varphi|_W=0,\ \ \varphi(a)=1.
\]
Let $\tau=I+u\otimes\varphi$. Since $u\ne0$, $\varphi(a)=1$, and $\varphi(u)=0$, we have $\tau\in T_n$, and
\begin{align*}
	\tau x&=x+\varphi(x)u=x,\ \ x\in U,\\
	\tau a&=a+\varphi(a)u=a+(a+b)=b.
\end{align*}
\end{proof}

\begin{remark}\label{r25}
Let $V$ be a vector space over $\mathbb F_q$, $U\subseteq V$ be a linear subspace, and  $a,b\in V\setminus U$ with $a\ne b$. A transvection $\tau$ satisfying $\tau|_U=\mathrm{id}_U$ and $\tau a=b$ exists if and only if
\[
a\notin U+\mathbb F_q(b-a).
\]
Indeed, if such a $\tau$ exists, then $(\tau-I)(b-a)=(\tau-I)^2a=0$, and therefore
\[
U+\mathbb F_q(b-a)\subseteq\ker(\tau-I),\ \ (\tau-I)a=b-a\ne0.
\]
Conversely, if $a\notin U+\mathbb F_q(b-a)$, choose $\varphi\in V^*$ such that
\[
\varphi|_{U+\mathbb F_q(b-a)}=0,\ \ \varphi(a)=1,
\]
then $\tau=I+(b-a)\otimes\varphi$ is the desired transvection.

In particular, when $q>2$, take $b=\lambda a$ with $\lambda\in\mathbb F_q\setminus\{0,1\}$. Then
\[
a=(\lambda-1)^{-1}(b-a)\in U+\mathbb F_q(b-a).
\]
Hence such a transvection does not exist, and Lemma~\ref{l24} cannot be directly generalized to arbitrary finite fields.
\end{remark}

The transvection constructed in Lemma \ref{l24} also yields a bijection between the corresponding two sets of completions. If $\tau|_U=\mathrm{id}_U$ and $\tau a=b$, then left multiplication by $\tau$ gives a bijection
\[
B(v_{<k},a)\longrightarrow B(v_{<k},b).
\]
Indeed, for
$g=(v_1,\ldots,v_{k-1},a,v_{k+1},\ldots,v_n)\in B(v_{<k},a)$, we have
\[
\tau g=(v_1,\ldots,v_{k-1},b,\tau v_{k+1},\ldots,\tau v_n)\in B(v_{<k},b),
\]
and invertibility of $\tau$ gives the inverse map.

\begin{lemma}\label{l26}
Fix $v_{<k}$, and take $a,b\in V\setminus U(v_{<k})$ with $a\ne b$. Let $\tau$ be a transvection obtained from Lemma \ref{l24}. Then
\[
\bigl(\Psi_{v_{<k}}(a)-\Psi_{v_{<k}}(b)\bigr)^2
\le
\frac1{|B(v_{<k},a)|}
\sum_{g\in B(v_{<k},a)}
\bigl(f(\tau g)-f(g)\bigr)^2.
\]
\end{lemma}

\begin{proof}
By the above bijection, we have
\begin{align*}
	\Psi_{v_{<k}}(a)&=\left(\frac1{|B(v_{<k},a)|}\sum_{g\in B(v_{<k},a)}f(g)^2\right)^{1/2},\\
	\Psi_{v_{<k}}(b)&=\left(\frac1{|B(v_{<k},a)|}\sum_{g\in B(v_{<k},a)}f(\tau g)^2\right)^{1/2}.
\end{align*}
In $L^2(B(v_{<k},a))$ equipped with the normalized counting measure, the reverse triangle inequality yields
\begin{align*}
	\bigl(\Psi_{v_{<k}}(a)-\Psi_{v_{<k}}(b)\bigr)^2
	&=\bigl(\|f\|_{L^2(B(v_{<k},a))}-\|f(\tau\,\cdot)\|_{L^2(B(v_{<k},a))}\bigr)^2\\
	&\le\|f(\tau\,\cdot)-f\|_{L^2(B(v_{<k},a))}^2\\
	&=\frac1{|B(v_{<k},a)|}\sum_{g\in B(v_{<k},a)}\bigl(f(\tau g)-f(g)\bigr)^2.
\end{align*}
\end{proof}

\begin{proposition}\label{p27}
For every real-valued $f$ on $G_n$, we have
\[
\Ent_{\pi_n}(f^2)
\le
\sum_{k=1}^n c_{m_k} Q_k(f).
\]
\end{proposition}

\begin{proof}
By Lemma~\ref{l24}, the map \(\Theta_k\) used in the definition of \(Q_k\) exists, and all the sets involved are finite. Fix \(v_{<k}\), By Lemma~\ref{l26}, we have
\begin{align*}
	D_k(v_{<k})
	&=\frac1{2m_k(m_k-1)}\sum_{\substack{a,b\in X(v_{<k})\\a\ne b}}\bigl(\Psi_{v_{<k}}(a)-\Psi_{v_{<k}}(b)\bigr)^2\\
	&\le\frac1{2m_k(m_k-1)}\sum_{\substack{a,b\in X(v_{<k})\\a\ne b}}
	\frac1{|B(v_{<k},a)|}\sum_{g\in B(v_{<k},a)}\bigl(f(\Theta_k(v_{<k},a,b)g)-f(g)\bigr)^2\\
	&=\frac12\sum_{(g,b)\in\Omega_k(v_{<k})}\bigl(f(\tau_k(g,b)g)-f(g)\bigr)^2\rho_k(g,b\mid v_{<k}),
\end{align*}
where $\Omega_k(v_{<k})$ denotes the part of $\Omega_k$ whose first $k-1$ columns are $v_{<k}$, and $\rho_k(\cdot\mid v_{<k})$ denotes the conditional probability measure of $\rho_k$ on this part. Taking the expectation with respect to $V_{<k}$, we obtain
\begin{align*}
	\mathbb E_{\pi_n}[D_k(V_{<k})]
	&\le\frac12\sum_{(g,b)\in\Omega_k}\bigl(f(\tau_k(g,b)g)-f(g)\bigr)^2\rho_k(g,b)\\
	&=Q_k(f).
\end{align*}
Combining this with Lemma~\ref{l23}, we obtain
\[
\Ent_{\pi_n}(f^2)\le\sum_{k=1}^n c_{m_k}\mathbb E_{\pi_n}[D_k(V_{<k})]\le\sum_{k=1}^n c_{m_k}Q_k(f).
\]
\end{proof}

\section{Averaging over left and right translations and proof of the main theorem}
The estimate in Section~\ref{s22} depends on the transvections chosen according to the values of the preceding columns. We now eliminate this dependence by averaging. Left translation makes each chosen transvection range over the entire conjugacy class $T_n$, while right translation makes the starting point of the random walk uniformly distributed over $G_n$. Therefore, the averaged local energy is exactly equal to the original Dirichlet form.

For $L,R\in G_n$, define
\[
f_{L,R}(g)=f(LgR),\ \ g\in G_n.
\]
Since left and right translations both preserve $\pi_n$, then
\[
\Ent_{\pi_n}(f_{L,R}^2)=\Ent_{\pi_n}(f^2).
\]

\begin{proposition}\label{p31}
For every $1\le k\le n$, we have
\[
\sum_{L\in G_n}\sum_{R\in G_n} Q_k(f_{L,R})\pi_n(L)\pi_n(R)
=\mathcal E_n(f,f).
\]
\end{proposition}

\begin{proof}
By definition, then
\[
f_{L,R}(\tau_k(g,b)g)-f_{L,R}(g)
=f(L\tau_k(g,b)gR)-f(LgR).
\]
Let
\[
S=L\tau_k(g,b)L^{-1},\ \ X=LgR,
\]
therefore
\[
f(L\tau_k(g,b)gR)-f(LgR)=f(SX)-f(X).
\]
We now show that, with respect to the probability mass function
\[
\pi_n(L)\pi_n(R)\rho_k(g,b), \ \ L,R\in G_n, \ \ (g,b)\in\Omega_k,
\]
the joint distribution of the random variables $(S,X)$ is $\mu_n\otimes\pi_n$.

First fix $(g,b)\in\Omega_k$, then $\tau_k(g,b)\in T_n$. Let
$\mathcal O(\tau_k)=\{A\tau_kA^{-1}:A\in G_n\}$ be its conjugacy class, and 
$C_{G_n}(\tau_k)=\{A\in G_n:A\tau_k=\tau_kA\}$ be its centralizer. Since $T_n$ is a conjugacy class, we have $\mathcal O(\tau_k)=T_n$. By the orbit--stabilizer theorem, then
\[
|\mathcal O(\tau_k)|=\frac{|G_n|}{|C_{G_n}(\tau_k)|}.
\]
For any $\sigma\in T_n$, the set $\{L\in G_n:L\tau_kL^{-1}=\sigma\}$
is a left coset of $C_{G_n}(\tau_k)$ and hence has cardinality $|C_{G_n}(\tau_k)|$. Therefore
\[
\sum_{\substack{L\in G_n:\, L\tau_kL^{-1}=\sigma}}\pi_n(L)
=\frac{|C_{G_n}(\tau_k)|}{|G_n|}
=\frac1{|\mathcal O(\tau_k)|}
=\frac1{|T_n|}
=\mu_n(\sigma).
\]
Hence $S=L\tau_k(g,b)L^{-1}$ has distribution $\mu_n$.

Next, fix $L,g,b$. Since the map $R\mapsto LgR$ is a bijection from $G_n$ to itself, when $R$ is distributed according to the uniform probability measure $\pi_n$, the distribution of $X=LgR$ is still $\pi_n$. To determine the joint distribution, take arbitrary functions $A:T_n\to\mathbb R$ and $B:G_n\to\mathbb R$. Then, we have
\begin{align*}
&\sum_{L\in G_n}\sum_{R\in G_n}\sum_{(g,b)\in\Omega_k}
A(S)B(X)\pi_n(L)\pi_n(R)\rho_k(g,b)\\
&=\sum_{L\in G_n}\sum_{(g,b)\in\Omega_k}
A(L\tau_k(g,b)L^{-1})
\left(\sum_{R\in G_n}B(LgR)\pi_n(R)\right)
\pi_n(L)\rho_k(g,b)\\
&=\left(\sum_{x\in G_n}B(x)\pi_n(x)\right)
\sum_{L\in G_n}\sum_{(g,b)\in\Omega_k}
A(L\tau_k(g,b)L^{-1})\pi_n(L)\rho_k(g,b)\\
&=\left(\sum_{x\in G_n}B(x)\pi_n(x)\right)
\left(\sum_{\sigma\in T_n}A(\sigma)\mu_n(\sigma)\right).
\end{align*}
Hence the joint distribution of $(S,X)$ is $\mu_n\otimes\pi_n$. Consequently,
\begin{align*}
\sum_{L\in G_n}\sum_{R\in G_n}Q_k(f_{L,R})\pi_n(L)\pi_n(R)
&=\frac12\sum_{\sigma\in T_n}\sum_{x\in G_n}
\bigl(f(\sigma x)-f(x)\bigr)^2\mu_n(\sigma)\pi_n(x)\\
&=\mathcal E_n(f,f).
\end{align*}
\end{proof}

We can now prove the main theorem of this paper.

\begin{proof}[Proof of Theorem \ref{t1}]
Applying Proposition~\ref{p27} to $f_{L,R}$ and averaging over $L,R\in G_n$, and using the invariance of entropy under left and right translations together with Proposition~\ref{p31}, we obtain
\begin{align*}
\Ent_{\pi_n}(f^2)
&=\sum_{L\in G_n}\sum_{R\in G_n}\Ent_{\pi_n}(f_{L,R}^2)\pi_n(L)\pi_n(R)\\
&\le \sum_{k=1}^n c_{m_k}
\sum_{L\in G_n}\sum_{R\in G_n}Q_k(f_{L,R})\pi_n(L)\pi_n(R)\\
&=\left(\sum_{k=1}^n c_{m_k}\right)\mathcal E_n(f,f).
\end{align*}
Therefore
\[
C_{\mathrm{LS}}(n)\le \sum_{k=1}^n c_{m_k}.
\]
Since $m_k=2^n-2^{k-1}\ge 2^{n-1}$ and
\[
c_m=\log m+O\left(\frac{\log m}{m}\right),
\]
we have
\begin{align*}
\sum_{k=1}^n c_{m_k}
&=\sum_{k=1}^n \log m_k
+O\left(\sum_{k=1}^n \frac{\log m_k}{m_k}\right)\\
&=\log |G_n|+O\left(\frac{n^2}{2^n}\right)\\
&=\log |G_n|+o(1).
\end{align*}
This proves the upper bound.

Next we prove the lower bound. Take the test function $f=\mathbf 1_{\{e\}}$, where $e=I$ is the identity element of $G_n$. Since $I\notin T_n$, we have $P_nf(e)=0$. Hence
\[
\mathcal E_n(f,f)=\langle f,(I-P_n)f\rangle_{\pi_n}
=\pi_n(e)(1-P_nf(e))
=\frac1{|G_n|}.
\]
On the other hand, $f^2=f$, and hence
\[
\Ent_{\pi_n}(f^2)=\frac1{|G_n|}\log |G_n|.
\]
It follows that
\[
C_{\mathrm{LS}}(n)\ge \log |G_n|.
\]
Combining the upper and lower bounds gives the first assertion of the theorem.

Finally, from the exact formula for the order of the group, we have
\[
|G_n|=|\GL_n(\mathbb F_2)|=|\SL_n(\mathbb F_2)|
=\prod_{j=0}^{n-1}(2^n-2^j)
=2^{n^2}\prod_{r=1}^n(1-2^{-r}).
\]
Then
\[
\log |G_n|
=n^2\log 2+\sum_{r=1}^n\log(1-2^{-r})
=n^2\log 2+O(1),
\]
by the infinite sum $\sum_{r\ge1}\log(1-2^{-r})$ converges. Therefore
\[
C_{\mathrm{LS}}(n)=n^2\log 2+O(1).
\]
\end{proof}

\section*{Acknowledgements}
The author thanks Professor Zipeng Wang and Professor Penghui Wang for helpful discussions, as well as their valuable comments on the manuscript.

\end{document}